\documentclass[oneside,english]{amsart}
\usepackage{mathtools, amssymb}
\usepackage[T1]{fontenc}
\usepackage[latin9]{inputenc}
\usepackage{geometry}
\usepackage{amsthm}

\usepackage{xcolor}

\theoremstyle{plain}
\newtheorem*{conjecture*}{\protect\conjecturename}
\theoremstyle{plain}
\newtheorem{cor}{\protect\corollaryname}
\theoremstyle{remark}
\newtheorem*{rem*}{\protect\remarkname}
\theoremstyle{plain}
\newtheorem{lem}{\protect\lemmaname}
\theoremstyle{plain}
\newtheorem{thm}{\protect\theoremname}
\theoremstyle{definition}
\newtheorem{defn}{\protect\definitionname}

\makeatother

\usepackage{babel}
\providecommand{\conjecturename}{Conjecture}
\providecommand{\corollaryname}{Corollary}
\providecommand{\definitionname}{Definition}
\providecommand{\lemmaname}{Lemma}
\providecommand{\remarkname}{Remark}
\providecommand{\theoremname}{Theorem}

\newtheorem*{example*}{\protect\examplename}

\providecommand{\examplename}{Example}

\makeatletter

\providecommand{\definitionname}{Definition}

\makeatother

\usepackage{babel}

\title{A Generalization of the Erd\H{o}s-Kac Theorem}

\begin{document}

\author {Matthew Levy \\ mlevy@skewfinancial.com \\  \\ Joseph Squillace \\ 
joseph.squillace@uvi.edu}

\begin{abstract}
Given a natural number $n$, let $\omega\left(n\right)$ denote the number
of distinct prime factors of $n$, let $Z$ denote a standard normal
variable, and let $P_{n}$ denote the uniform distribution on $\left\{ 1,\ldots,n\right\} $.
The Erd\H{o}s-Kac Theorem states that if $N\left(n\right)$ is a uniformly distributed variable on $\lbrace 1,\ldots,n \rbrace$, then $\omega\left(N\left(n\right)\right)$ is asymptotically normally distributed as $n\to \infty$ with both mean and variance equal to $\log \log n$. The contribution of this paper is a generalization of the Erd\H{o}s-Kac Theorem to a larger class
of random variables by considering perturbations of the uniform probability
mass $1/n$ in the following sense. Denote by $\mathbb{P}_{n}$
a probability distribution on $\left\{ 1,\ldots,n\right\} $ given
by $\mathbb{P}_{n}\left(i\right)=1/n+\varepsilon_{i,n}$.
We provide sufficient conditions on $\varepsilon_{i,n}$ so that the number of distinct prime factors of a $\mathbb{P}_{n}$-distributed random variable is asymptotically normally distributed, as $n\to \infty$, with both mean and variance equal to $\log \log n$. Our main result is applied to prove that the number of distinct
prime factors of a positive integer with the Harmonic$\left(n\right)$
distribution also tends to the normal distribution, as $n\to \infty$. In addition, we explore sequences of distributions on the natural numbers such that $\omega(n)$ is normally distributed in the limit. In addition, one of our theorems and its corollaries generalize a result from the literature involving the limit of $Zeta\left(s\right)$ distributions as the parameter $s \to 1$.
\end{abstract}

\maketitle

\section{Introduction}
Given a natural number $n$, the number of distinct prime factors
of $n$ is denoted $\omega\left(n\right)$. The function $\omega$ may be written as 
$\omega\left(n\right)=\sum_{p\vert n}1$, 
where the sum is over all prime factors of $n$. In 1917, Hardy and Ramanujan (p. 270 of \cite{key-5}) proved
that the number of distinct prime factors of a natural number $n$
is about $\log\log n$. In particular, they showed that the normal
order of $\omega\left(n\right)$ is $\log\log n$; i.e., for every
$\varepsilon>0$, the proportion of the natural numbers for which
the inequalities 
\[
\left(1-\varepsilon\right)\log\log n\le\omega\left(n\right)\le\left(1+\varepsilon\right)\log\log n
\]
do not hold tends to $0$ as $n\to\infty$--with a typical error of size $\sqrt{\log \log n}$. Informally speaking, the
Erd\H{o}s-Kac Theorem generalizes\footnote{While the Hardy-Ramanujan theorem provides information about the average behavior of the number of prime factors of a natural number, the Erd\H{o}s-Kac theorem offers a more detailed probabilistic description of their distribution, by taking into account not just the average number but also the variability around that average.} the Hardy-Ramanujan Theorem by showing that $\omega\left(n\right)$ is approximately distributed as
\[\log\log n+Z\sqrt{\log\log n}\] 
for large $n$, where $Z$ denotes a standard normal variable. 
More precisely, the Erd\H{o}s-Kac Theorem is the following result (p. 738 of \cite{key-4}).

\begin{thm}\label{thm:1}
Let $n > 1$. Let $P_{n}$ denote the uniform distribution on $\left\{ 1,2,,\ldots,n\right\}$, 
and let $Z$ denote a standard normal variable. As $n\to\infty$, 
\[
P_{n}\left(m\le n:\omega\left(m\right)-\log\log n\le x\left(\log\log n\right)^{1/2}\right)\to\mathbb{P}\left(Z\le x\right).
\]
\end{thm}

The contribution of this paper is to extend the Erd\H{o}s-Kac Theorem to
a larger class of random variables on the set $\left[n\right]\coloneqq\left\{ 1,2,\ldots,n\right\}$
which also have, asymptotically, $\log\log n+Z\cdot\sqrt{\log\log n}$ many distinct prime factors.

\subsection{A Generalization of Erd\H{o}s-Kac Theorem for $\omega\left(\cdot \right)$}

Define a probability mass function (PMF) $\mathbb{P}_{n}$
on $\left[n\right]$ given by 
\begin{equation}
\mathbb{P}_{n}\left(i\right)=\frac{1}{n}+\varepsilon_{i,n}.\label{eq:1}
\end{equation}
Due to the axioms of probability, the terms $\varepsilon_{i,n}, 1\le i \le n$, satisfy the constraints
\begin{equation}
\sum_{i=1}^{n}\varepsilon_{i,n}= 0; \,\, n \in \mathbb{N}\label{eq:2}
\end{equation}
and
\begin{equation}
\varepsilon_{i,n}\in \left[-\frac{1}{n},1-\frac{1}{n}\right]; \,\, n \in \mathbb{N}, 1 \le i\le n.\label{eq:3}
\end{equation}

The motivation for defining $\mathbb{P}_{n}$ in terms of the uniform
distribution is due to Durrett's proof (Theorem 3.4.16 in \cite{key-3})
of the Erd\H{o}s-Kac Theorem. Replacing the uniform distribution
$P_{n}$ with the distribution $\mathbb{P}_{n}$ in Durrett's
proof naturally yields some constraints that the terms $\varepsilon_{i,n},1\le i\le n$,
must satisfy in order to conclude that an integer-valued $\mathbb{P}_{n}$-distributed random variable 
has about $\log\log n+Z\sqrt{\log\log n}$
distinct prime factors. Our main result is the following theorem,
where $\left\lfloor \cdot\right\rfloor $ denotes the floor function.

\begin{thm}{(Generalized Erd\H{o}s-Kac Theorem for $\omega$)}\label{thm:2}
Let $Z$ denote a standard normal variable. Suppose the following statements are true.

\begin{itemize}
\item There exists a constant $C\in\mathbb{R}$ such that for all $n > 1$ and for all primes
$p$ with $p>n^{1/\log\log n}$,
\begin{equation}
\sum_{l=1}^{\left\lfloor n/p\right\rfloor}\varepsilon_{lp,n}\le\frac{C}{p}.\label{eq:4}
\end{equation}
\item There exists a constant $D\in\mathbb{R}$ such that
\begin{equation}
\sum_{l=1}^{\left\lfloor \frac{n}{p_{1}\cdots p_{k}}\right\rfloor }\varepsilon_{lp_{1}\cdots p_{k},n}\le\frac{D}{n}\label{eq:5}
\end{equation}
for all $n > 1$, for each $k$, and for all $k$-tuples $\left(p_{1},\ldots,p_{k}\right)$
consisting of distinct primes of size at most $n^{1/\log\log n}$. 
\item For any prime $p$, 
\begin{equation}
\lim_{n \to \infty} \sum_{l=1}^{\lfloor n/p \rfloor} \varepsilon_{lp, n} = 0.\label{eq:6}
\end{equation}

\end{itemize}
Let $\mathbb{P}_{n}^{*}$ denote the PMF obtained
by imposing the constraints $\left(\ref{eq:4}-\ref{eq:6}\right)$ on the PMF$\,$ $\mathbb{P}_{n}$
given by $\mathbb{P}_{n}\left(i\right)=\frac{1}{n}+\varepsilon_{i,n}$. As $n\to\infty,$ 
\begin{align*}
\mathbb{P}_{n}^{*}\left(m\le n:\omega\left(m\right)-\log\log n\le x\left(\log\log n\right)^{1/2}\right) & \to\mathbb{P}\left(Z\le x\right).
\end{align*}
\end{thm}

\begin{rem*}
If $\varepsilon_{i,n}=0$ for all $i\le n$, then $\mathbb{P}_{n}^{*}=P_{n}$ and Theorem \ref{thm:1} is obtained.
\end{rem*}

\subsection{Outline}
The proof of Theorem \ref{thm:2} is provided in $\S$2; the proof applies the method of moments and is motivated by Durrett's proof of the Erd\H{o}s-Kac Theorem (Theorem 3.4.16 in
\cite{key-3}). Moreover, in $\S$2, the constraints $\left(\ref{eq:4}-\ref{eq:6}\right)$ are applied to ensure that $\mathbb{P}_n^{*}$ also satisfies Durrett's method of moment bounds. In $\S\S$2.1, Theorem \ref{thm:2} is applied to show that the number of distinct prime factors of a random natural number chosen according to the Harmonic$\left(n\right)$ distribution is asymptotically normally distributed with both mean and variance equal to $\log \log n$. 
In $\S$3, Theorem \ref{thm:2} is used to prove statements about convex combinations of distributions satisfying constraints $\left(\ref{eq:4}-\ref{eq:6}\right)$. In $\S$4, we define conditions that ensure $\omega(X_{j}\left(n\right))$ is asymptotically normally distributed, with mean and variance both equal to $\log \log n$, for a sequence of random variables $\left(X_{j}\left(n\right)\right)_{j\ge 1}$ as $j\to \infty$ and $n\to \infty$.
In $\S$5, Theorem \ref{thm:2} is applied to show that the number of distinct prime factors of a randomly chosen integer according to any of the following distributions has the same limiting distribution as the case of a uniform variable:
\begin{itemize}
\item Any convex combination of the Harmonic$\left(n\right)$ and uniform$\left(n\right)$ distributions,
\item The $Zeta\left(s\right)$ distribution as $s\to 1$,
\item The Logarithmic$\left(s\right)$ distribution as $s\to 1$,
\item A geometric power series distribution as $s \to 1$,
\item A Logarithmic-Zeta$\left(s,\alpha\right)$ distribution as $(s,\alpha) \to (1,1)$.
\end{itemize}

\section{Proving Theorem \ref{thm:2}}

Define $\alpha_{n}\coloneqq n^{1/\log\log n}$.
\begin{lem}\label{lem:1}
As $n\to \infty$
\[
\left(\sum_{\alpha_{n}<p\le n}\left(\frac{1}{p}+\sum_{l=1}^{\left\lfloor n/p \right\rfloor }\varepsilon_{lp,n}\right)\right)/\left(\log\log n\right)^{1/2}\to0.
\]
\end{lem}

\begin{proof}
Given $n$ and any prime $p$ with $p>\alpha_{n}$, we have
\[
-\frac{\left\lfloor n/p\right\rfloor }{n}\overset{\left(\ref{eq:3}\right)}{\le}\sum_{l=1}^{\left\lfloor n/p\right\rfloor }\varepsilon_{lp,n}\overset{\left(\ref{eq:4}\right)}{\le}\frac{C}{p}.
\]
 Therefore,
\begin{equation}
\frac{1}{p}+\sum_{l=1}^{\left\lfloor n/p\right\rfloor }\varepsilon_{lp,n}\in\left[0,\frac{C+1}{p}\right]\label{eq:7}
\end{equation}
for all $n$. Thus,
\[
\left(\sum_{\alpha_{n}<p\le n}\left(\frac{1}{p}+\sum_{l=1}^{\left\lfloor n/p\right\rfloor }\varepsilon_{lp,n}\right)\right)/\left(\log\log n\right)^{1/2}\to0
\]
due to $\left(\ref{eq:7}\right)$ along with the fact that Durrett (p.135 of \cite{key-3}) shows 
\[
\left(\sum_{\alpha_{n}<p\le n}\frac{1}{p}\right)/\left(\log\log n\right)^{1/2}\to0.
\]
\end{proof}
The following lemma is proved by Durrett (p. 156 of \cite{key-3}).
\begin{lem}
If $\varepsilon>0$, then $\alpha_{n}\le n^{\varepsilon}$ for large
$n$ and hence 
\begin{equation}
\frac{\alpha_{n}^{r}}{n}\to0\label{eq:8}
\end{equation}
for all $r<\infty$.
\end{lem}

\begin{proof}[Proof of Theorem \ref{thm:2}]
Given a natural number $m$ and a prime $p$, define $\delta_{p}\left(m\right)=1$ if $p$ divides $m$, and $0$ otherwise. Let
\[
g_{n}\left(m\right)=\sum_{p\le\alpha_{n}}\delta_{p}\left(m\right)
\]
denote the number of distinct prime factors of $m$ of size at most $\alpha_n$,
and let $\mathbb{E}_{n}$ denote expectation with respect to $\mathbb{P}_{n}^{*}$.
Then

\begin{align*}
\mathbb{E}_{n}\left(\sum_{\alpha_{n}<p\le n}\delta_{p}\right) & =\sum_{\alpha_{n}<p\le n}\mathbb{P}_{n}^{*}\left(m:\delta_{p}\left(m\right)=1\right)\\
 & =\sum_{\alpha_{n}<p\le n}\sum_{l=1}^{\left\lfloor n/p\right\rfloor }\mathbb{P}_{n}^{*}\left(m:m=lp\right)\\
 & \overset{\left(\ref{eq:1}\right)}{=}\sum_{\alpha_{n}<p\le n}\sum_{l=1}^{\left\lfloor n/p\right\rfloor }\left(\frac{1}{n}+\varepsilon_{lp,n}\right)\\
 & \le\sum_{\alpha_{n}<p\le n}\left(\frac{1}{p}+\sum_{l=1}^{\left\lfloor n/p\right\rfloor }\varepsilon_{lp,n}\right),
\end{align*}
so by Lemma \ref{lem:1} it suffices to prove Theorem \ref{thm:2} for $g_{n}$; i.e.,
replacing $\omega\left(m\right)$ with $g_{n}\left(m\right)$
in the statement of Theorem 2 does not affect the limiting distribution. 

Consider a sequence $\left( X_p \right)_{p\ge 2}$ of independent Bernoulli random variables with prime-valued indices such that $\mathbb{P}\left(X_p = 1\right) = 1/p$ and $\mathbb{P}\left(X_p = 0 \right) = 1-1/p$. Note that 
\[\mathbb{E}\left(\delta_p\right)=\frac{\left\lfloor n/p \right\rfloor}{n}+\sum_{l=1}^{\left\lfloor n/p \right\rfloor}\varepsilon_{lp, n} \overset{\left(6\right)}{\to} 1/p\]
as $n \to \infty$.
Let
\begin{eqnarray*}
S_{n} & \coloneqq & \sum_{p\le\alpha_{n}}X_{p},\\
b_{n} & \coloneqq & \mathbb{E}\left(S_{n}\right),\\
a_{n}^{2} & \coloneqq & \text{Var}\left(S_{n}\right).
\end{eqnarray*}
By Lemma \ref{lem:1}, $b_{n}$ and $a_{n}^{2}$ are both $\log\log n+o\left(\left(\log\log n\right)^{1/2}\right)$,
so it suffices to show
\[\mathbb{P}_{n}^{*}\left(m:g_{n}\left(m\right)-b_{n}\le xa_{n}\right)\to\mathbb{P}\left(Z\le x\right).\]
An application of Theorem 3.4.10 of \cite{key-3} shows 
\[\left(S_{n}-b_{n}\right)/a_{n}\to Z,\]
and since $\left|X_{p}\right|\le1$, it follows from Durrett's second
proof of Theorem 3.4.10 \cite{key-3} that 
\[\mathbb{E}\left(\left(S_{n}-b_{n}\right)/a_{n}\right)^{r}\to\mathbb{E}\left(Z^{r}\right)\]
for all $r$. Using the notation from that proof (and replacing $i_{j}$
by $p_{j}$) it follows that
\begin{align*}
\mathbb{E}\left(S_{n}^{r}\right) & =\sum_{k=1}^{r}\sum_{r_{i}}\frac{r!}{r_{1}!\cdots r_{k}!}\frac{1}{k!}\sum_{p_{j}}\mathbb{E}\left(X_{p_{1}}^{r_{1}}\cdots X_{p_{k}}^{r_{k}}\right),
\end{align*}
where the sum $\sum_{r_{i}}$ extends over all $k$-tuples of positive
integers for which $r_{1}+\cdots+r_{k}=r$, and $\sum_{p_{j}}$ extends
over all $k$-tuples of distinct primes in $\left[n\right]$. Since
$X_{p}\in\left\{ 0,1\right\} $, the summand in $\sum_{p_{j}}\mathbb{E}\left(X_{p_{1}}^{r_{1}}\cdots X_{p_{k}}^{r_{k}}\right)$
is
\[
\mathbb{E}\left(X_{p_{1}}\cdots X_{p_{k}}\right)=\frac{1}{p_{1}\cdots p_{k}}
\]
by independence of the $X_{p}$'s. Moreover,
\begin{align*}
\mathbb{E}_{n}\left(\delta_{p_{1}}\cdots\delta_{p_{k}}\right) & \le\mathbb{P}_{n}\left(m:\delta_{p_{1}}\left(m\right)=\delta_{p_{2}}\left(m\right)\cdots=\delta_{p_{k}}\left(m\right)=1\right)\\
 & =\sum_{l=1}^{\left\lfloor \frac{n}{p_{1}\cdots p_{k}}\right\rfloor }\mathbb{P}_{n}\left(m:m=lp_{1}\cdots p_{k}\right)\\
 & \overset{\left(1\right)}{=}\sum_{l=1}^{\left\lfloor \frac{n}{p_{1}\cdots p_{k}}\right\rfloor }\left(\frac{1}{n}+\varepsilon_{lp_{1}\cdots p_{k},n}\right)\\
 & =\frac{\left\lfloor \frac{n}{p_{1}\cdots p_{k}}\right\rfloor }{n}+\sum_{l=1}^{\left\lfloor \frac{n}{p_{1}\cdots p_{k}}\right\rfloor }\varepsilon_{lp_{1}\cdots p_{k},n}\\
 & \overset{\left(5\right)}{\le}\frac{\left\lfloor \frac{n}{p_{1}\cdots p_{k}}\right\rfloor }{n}+\frac{D}{n}.\\
\end{align*}
The two terms $\mathbb{E}\left(X_{p_{1}}\cdots X_{p_{k}}\right)$ and $\mathbb{E}_{n}\left(\delta_{p_{1}}\cdots\delta_{p_{k}}\right)$ differ by at most 
\[
\text{max}\left\{ \frac{1}{p_{1}\cdots p_{k}}-\frac{\left\lfloor \frac{n}{p_{1}\cdots p_{k}}\right\rfloor }{n}-\frac{D}{n},\frac{\left\lfloor \frac{n}{p_{1}\cdots p_{k}}\right\rfloor }{n}+\frac{D}{n}-\frac{1}{p_{1}\cdots p_{k}}\right\} \le \text{max}\left\{\frac{1-D}{n},\frac{D}{n}\right\}.
\]
Therefore, the two $r$th moments differ by
{\small{}
\begin{eqnarray*}
\left|\mathbb{E}\left(S_{n}^{r}\right)-\mathbb{E}_{n}\left(g_{n}^{r}\right)\right| 
 & \le & \sum_{k=1}^{r}\sum_{r_{i}}\frac{r!}{r_{1}!\cdots r_{k}!}\frac{1}{k!}\sum_{p_{j}}\max \left\{\frac{1-D}{n},\frac{D}{n}\right\}.\\
 & \le &\max \left\{\frac{1-D}{n},\frac{D}{n}\right\}\cdot\left(\sum_{p\le\alpha_{n}}1\right)^{r}\\
 & \le & \max \left\{\frac{1-D}{n},\frac{D}{n}\right\}\cdot \alpha_{n}^{r}\\
 & \overset{\left(\ref{eq:8}\right)}{\to} & 0.
\end{eqnarray*}
}{\small\par}
\noindent Using binomial expansions and the inequality above, we see that
\begin{align*}
\left|\mathbb{E}\left( \left( \left(S_{n} - b_{n} \right) / a_n\right)^{r}\right)-\mathbb{E}_{n}\left( \left( \left(g_{n} - b_{n} \right) / a_{n}\right)^{r}\right) \right| & = 
\left|1/a_{n}^{r}\right|\left|\mathbb{E}\left( \left(S_{n} - b_{n} \right)^{r}\right)-\mathbb{E}_{n}\left( \left( g_{n} - b_{n} \right)^{r}\right) \right|\\
&  
\leq \left| 1/a_{n}^{r}  \right| \cdot \max \left\{\frac{1-D}{n},\frac{D}{n}\right\} \sum_{k=0}^{r} \binom{r}{k} \alpha_{n}^{k} b_{n}^{r-k} \\
& = \left| 1/a_{n}^{r}  \right| \cdot \max \left\{\frac{1-D}{n},\frac{D}{n}\right\} \left( \alpha_{n} + b_{n} \right)^{r}.
\end{align*}
\noindent Therefore, since $b_{n} \leq \alpha_{n}$, we have
\[
\left|\mathbb{E}\left( \left( \left(S_{n} - b_{n} \right) / a_n\right)^{r}\right)-\mathbb{E}_{n}\left( \left( \left(g_{n} - b_{n} \right) / a_{n}\right)^{r}\right) \right| \overset{\left(\ref{eq:8}\right)}{\to} 0
\]
for all $r$ as well.
\noindent Since $\mathbb{E}\left( \left( \left(S_{n} - b_{n} \right) / a_{n}\right)^{r}\right) \to \mathbb{E}\left(Z^{r}\right)$ for all $r$, this completes the proof of Theorem \ref{thm:2}.
\end{proof}

The following definition is based on Theorem \ref{thm:2}.
\begin{defn}\label{def:1}
We refer to distributions satisfying constraints $\left(\ref{eq:4}-\ref{eq:6}\right)$ as \textbf{E-K distributions}.
\end{defn}

\subsection{The Harmonic Distribution}
Now we will apply Theorem \ref{thm:2} to show that the harmonic distributions are E-K distributions. Given $n\in\mathbb{N},$ consider an integer in $\left[n\right]$ chosen according to the Harmonic$\left(n\right)$ distribution, whose PMF is given by
\[ H_{n}\left(i\right)\coloneqq \frac{1}{i\sum_{j=1}^{n}\frac{1}{j}},  1\le i \le n.
\]
If $i\in\left[n\right]$, then equation $\left(\ref{eq:1}\right)$ implies
\[\varepsilon_{i,n}=\frac{1}{i\sum_{j=1}^{n}\frac{1}{j}}-\frac{1}{n}.\]
Therefore,
\begin{align*}
\sum_{l=1}^{\left\lfloor n/p \right\rfloor }\varepsilon_{lp,n} & =\sum_{l=1}^{\left\lfloor n/p \right\rfloor }\left(\frac{1}{lp\sum_{i=1}^{n}\frac{1}{i}}-\frac{1}{n}\right)\\
 & =\frac{\sum_{i=1}^{\left\lfloor n/p \right\rfloor }\frac{1}{i}}{p\sum_{l=1}^{n}\frac{1}{i}}-\frac{\lfloor n/p
 \rfloor }{n}\\
 & \le\frac{1}{p}-\frac{\left\lfloor
 n/p
 \right\rfloor }{n}\\
 & \le\frac{1}{p},
\end{align*}
so $\left(\ref{eq:4}\right)$ holds with $C=1$. Moreover,
\begin{align*}
\sum_{l=1}^{\left\lfloor \frac{n}{p_{1}\cdots p_{k}}\right\rfloor }\varepsilon_{lp_{1}\cdots p_{k},n} & =\sum_{l=1}^{\left\lfloor \frac{n}{p_{1}\cdots p_{k}}\right\rfloor }\left(\frac{1}{lp_{1}\cdots p_{k}\sum_{i=1}^{n}\frac{1}{i}}-\frac{1}{n}\right)\\
 & =\frac{\sum_{i=1}^{\left\lfloor \frac{n}{p_{1}\cdots p_{k}}\right\rfloor }\frac{1}{i}}{p_{1}\cdots p_{k}\sum_{i=1}^{n}\frac{1}{i}}-\frac{\left\lfloor \frac{n}{p_{1}\cdots p_{k}}\right\rfloor }{n}\\
 & \le\frac{1}{p_{1}\cdots p_{k}}-\left(\frac{\frac{n}{p_{1}\cdots p_{k}}-1}{n}\right)\\
 & =1/n,
\end{align*}
so $\left(\ref{eq:5}\right)$ holds with $D=1$. Finally,
\begin{align*}
\lim_{n \to \infty} \sum_{l=1}^{\left\lfloor n/p
\right\rfloor }\varepsilon_{lp,n} & = \lim_{n \to \infty} \left(\frac{\sum_{i=1}^{\left\lfloor n/p
\right\rfloor }\frac{1}{i}}{p\sum_{l=1}^{n}\frac{1}{i}}-\frac{\lfloor 
n/p \rfloor }{n}\right)\\
& \sim \frac{1}{p}\frac{\log \left(n/p\right)}{\log n} - \frac{1}{p}\\
 & \to  1/p - 1/p\\
 & =  0,
\end{align*}
so constraint 
$\left( \ref{eq:6}
\right)$ holds. By Theorem \ref{thm:2}, this shows that the number of distinct prime factors of an $H_{n}$-distributed random variable is asymptotically normally distributed, as $n \to \infty$, with both mean and variance equal to 
$\log \log n$.

\section{Convex Combinations of E-K Distributions}
The following theorem shows that any convex sum of two E-K distributions is also an E-K distribution.

\begin{thm}{(Convex Combinations for Erd\H{o}s-Kac (CLT))}\label{thm:3}
Let $n > 1$ and $0\le \lambda \le 1$ be fixed. Suppose $d_{1,n}$ and $d_{2,n}$ are two PMFs on $\left[n\right]$ satisfying the constraints $\left(\ref{eq:4}-\ref{eq:6}\right)$ for all $n > 1$.
Then any PMF of the form
\begin{equation}\label{eq:9}
\mathbb{P}^{*}_{n}\left(i\right) = \lambda d_{1,n}\left(i\right) + \left(1-\lambda \right)d_{2,n}\left(i\right), 1\le i \le n
\end{equation}
also satisfies constraints $\left(\ref{eq:4}-\ref{eq:6}\right)$.
In particular,
\begin{align*}
\mathbb{P}_{n}^{*}\left(m\le n:\omega\left(m\right)-\log\log n\le x\left(\log\log n\right)^{1/2}\right) & \to\mathbb{P}\left(Z\le x\right)
\end{align*}
as $n \to \infty$.
\end{thm}

\begin{proof}
Denote
\[
\varepsilon_{i,n}=d_{1,n}\left(i\right)-1/n,
\]
\[
\varepsilon'_{i,n}=d_{2,n}\left(i\right)-1/n,
\]
and 
\[
\hat{\varepsilon}_{i,n}=\mathbb{P}^{*}_{n}\left(i\right)-1/n.
\]
We have 
\begin{align*}
\sum_{l=1}^{\left\lfloor n/p\right\rfloor }\hat{\varepsilon}_{lp,n} & =\sum_{l=1}^{\left\lfloor n/p \right\rfloor }\left(\lambda d_{1,n}\left(lp\right)+\left(1-\lambda \right)d_{2,n}\left(lp\right)-\frac{1}{n}\right)\\
 & \overset{\left(\ref{eq:1}\right)}{=} \sum_{l=1}^{\left\lfloor n/p \right\rfloor }\left(\lambda\left(\varepsilon_{lp,n} +1/n\right)+\left(1-\lambda \right)\left(\varepsilon'_{lp,n}+1/n\right)-\frac{1}{n}\right)\\
 & = \sum_{l=1}^{\lfloor n/p\rfloor}\left(\lambda\varepsilon_{lp,n}+\left(1-\lambda\right)\varepsilon_{lp,n}' \right) \\
 & =\lambda\sum_{l=1}^{\left\lfloor n/p \right\rfloor }\varepsilon_{lp,n}+\left(1-\lambda\right)\sum_{l=1}^{\left\lfloor n/p\right\rfloor }\varepsilon'_{lp,n}\\
 & \le \frac{\lambda C_{1}+\left(1-\lambda\right)C_{2}}{p},
\end{align*}
where the latest inequality is obtained by applying constraint
$\left(\ref{eq:4}\right)$ to both $d_{1,n}$ and $d_{2,n}$; therefore, 
$\left(\ref{eq:4}\right)$ holds for the PMF given by $\left(\ref{eq:9}\right)$ with $C=\lambda C_{1} + \left(1-\lambda\right)C_{2}$. Similarly, 
\begin{align*}
\sum_{l=1}^{\left\lfloor \frac{n}{p_{1}\cdots p_{k}}\right\rfloor }\hat{\varepsilon}_{lp_{1}\cdots p_{k},n} & =\sum_{l=1}^{\left\lfloor \frac{n}{p_{1}\cdots p_{k}}\right\rfloor }\left(\lambda d_{1,n}\left(lp_{1}\cdots p_{k}\right)+\left(1-\lambda \right)d_{2,n}\left(lp_{1}\cdots p_{k}\right)-\frac{1}{n}\right)\\
 & \overset{\left(\ref{eq:1}\right)}{=} \sum_{l=1}^{\left\lfloor \frac{n}{p_{1}\cdots p_{k}}\right\rfloor }\left(\lambda \left(\varepsilon_{lp_{1}\cdots p_{k},n} +1/n\right)+\left(1-\lambda \right)\left(\varepsilon'_{lp_{1}\cdots p_{k},n}+1/n\right)-\frac{1}{n}\right)\\
 & = \sum_{l=1}^{\lfloor \frac{n}{p_{1}\cdots p_{k}}\rfloor}\left(\lambda\varepsilon_{lp_{1}\cdots p_{k},n}+\left(1-\lambda \right)\varepsilon_{lp_{1}\cdots p_{k},n}' \right) \\
 & = \lambda \sum_{l=1}^{\left\lfloor \frac{n}{p_{1}\cdots p_{k}}\right\rfloor }\varepsilon_{lp_{1}\cdots p_{k},n}+\left(1-\lambda\right)\sum_{l=1}^{\left\lfloor \frac{n}{p_{1}\cdots p_{k}}\right\rfloor }\varepsilon'_{lp_{1}\cdots p_{k},n}\\
 &\le \frac{\lambda D_{1} + \left(1-\lambda\right)D_{2}}{n},
\end{align*}
where the latest inequality is obtained by applying constraint
$\left(\ref{eq:5}\right)$ to both $d_{1,n}$ and $d_{2,n}$; therefore, 
$\left(\ref{eq:5}\right)$ holds for the PMF given by $\left(\ref{eq:9}\right)$ with $D=\lambda D_{1} + \left(1 - \lambda\right) D_{2}$. Furthermore, 
\begin{align*}
    \lim_{n \to \infty}\sum_{l=1}^{\left\lfloor n/p \right\rfloor }\hat{\varepsilon}_{lp,n} 
& = \lambda \lim_{n \to \infty}\sum_{l=1}^{\lfloor n/p\rfloor}\varepsilon_{lp,n}+\left(1-\lambda\right)\lim_{n\to \infty}\sum_{l=1}^{\left\lfloor n/p\right\rfloor}\varepsilon_{lp,n}'\\
& = \lambda\cdot 0 + \left(1-\lambda\right)\cdot 0 \\
& = 0,
\end{align*}
where the second equation uses the fact that the distributions $d_1$ and $d_2$ satisfy constraint $\left(\ref{eq:6}\right)$.
Therefore, $\left(\ref{eq:6}\right)$ holds for the PMF given by $\left(\ref{eq:9}\right)$.

\end{proof}
By Theorem \ref{thm:3} and mathematical induction, we obtain the following.
\begin{cor}
Suppose $d_{1,n}, \ldots, d_{k,n}$ are PMFs on $\left[n\right]$ satisfying the constraints $\left(\ref{eq:4}-\ref{eq:6}\right)$,
and suppose $\left(\lambda _{j} \right)_{j=1}^{k}$ is a sequence of nonnegative numbers which add to unity. The PMF
\[
\mathbb{P}^{*}_{n}\left(i\right) = \sum_{j=1}^{k}\lambda_{j} d_{j,n}(i)
\]
also satisfies constraints $\left(\ref{eq:4}-\ref{eq:6}\right)$. Therefore, 
\begin{align*}
\mathbb{P}_{n}^{*}\left(m\le n:\omega\left(m\right)-\log\log n\le x\left(\log\log n\right)^{1/2}\right) & \to\mathbb{P}\left(Z\le x\right)
\end{align*}
as $n \to \infty$.
\end{cor}

We now present an alternate proof to the remark under Corollary 3.4 in \cite{key-1}; their proof uses both the Tauberian theorems for the harmonic distribution and the uniform distribution simultaneously. We present a proof based on Theorems \ref{thm:2} and \ref{thm:3} of this paper.

\begin{cor}{(Convex Combinations of Harmonic Distribution and Uniform Distribution)}\label{cor:2}
Let $\lambda \in \left[0,1\right]$ and consider the PMF given by
\[
\mathbb{P}^{*}_{n}\left(i\right) = \frac{\lambda}{ih_{n}} + \frac{1-\lambda}{n}; 1\le i \le n,
\]
where $h_{n}$ is the $n^{th}$ harmonic number given by $h_{n} \coloneq \sum_{i=1}^{n}1/i$, then
\[
\mathbb{P}_{n}^{*}\left(m\le n:\omega\left(m\right)-\log\log n\le x\left(\log\log n\right)^{1/2}\right) \to \mathbb{P}\left(Z\le x\right)
\]
as $n \to \infty$.
\end{cor}

\begin{proof}
As noted in the remark after Theorem \ref{thm:2}, the uniform distribution satisfies constraints $\left(\ref{eq:4}-\ref{eq:6}\right)$, and in $\S\S$2.1 we showed the harmonic distribution satisfies constraints $\left(\ref{eq:4}-\ref{eq:6}\right)$. By Theorem $\ref{thm:3}$, their convex sums also satisfy constraints $\left(\ref{eq:4}-\ref{eq:6}\right)$.
\end{proof}

The following theorem shows that if given the PMF $1/n+\varepsilon_{i,n}$ of a sequence of E-K distributions, under additional constraints, the PMF $1/n-\varepsilon_{i,n}$ is also the PMF of a sequence of E-K distributions.

\begin{thm}{(Reflection Theorem)}
Let $\mathbb{P}_{n}$ satisfy constraints $\left(\ref{eq:4}\right)$ and $\left(\ref{eq:6}\right)$ for all $n > 1$. In addition, assume $\varepsilon_{i,n}\in \left[-\frac{1}{n}, \frac{1}{n}\right]$ for all $i$ and for all $n$, and that for all $k$-tuples of distinct primes $p_{1}, \ldots, p_{k}$, there exists a constant $D\ge 0$ such that
\begin{equation}
-\frac{D}{n} \le \sum_{l=1}^{\lfloor \frac{n}{p_{1}\cdots p_{k}}\rfloor}\varepsilon_{lp_{1}\cdots p_{k},n} \le \frac{D}{n}. \label{eq:10}
\end{equation}
If
\[
\mathbb{P}_{n}^{*}\left(i\right) \coloneqq \frac{1}{n}-\varepsilon_{i,n},
\]
then
\begin{align*}
\mathbb{P}_{n}^{*}\left(m\le n:\omega\left(m\right)-\log\log n\le x\left(\log\log n\right)^{1/2}\right) & \to\mathbb{P}\left(Z\le x\right)
\end{align*}
as $n \to \infty$.
\end{thm}

\begin{proof}
Define $\varepsilon_{i,n}' \coloneq -\varepsilon_{i,n}$, so
\begin{align*}
\sum_{l=1}^{\left\lfloor n/p\right\rfloor }\varepsilon_{lp,n}' & =-\sum_{l=1}^{\left\lfloor n/p \right\rfloor}\varepsilon_{i,n} \\
& \overset{\left(\ref{eq:3}\right)}{\le}\frac{\lfloor n/p \rfloor}{n}\\
& \le \frac{1}{p}
\end{align*}
so constraint $\left(\ref{eq:4}\right)$ holds with $C=1$.
Similarly,
\begin{align*}
\sum_{l=1}^{\left\lfloor \frac{n}{p_{1}\cdots p_{k}}\right\rfloor }\varepsilon_{lp_{1}\cdots p_{k},n}' & =-\sum_{l=1}^{\left\lfloor \frac{n}{p_{1}\cdots p_{k}} \right\rfloor}\varepsilon_{lp_{1}\cdots p_{k},n} \overset{\left(\ref{eq:10}\right)}{\leq} \frac{D}{n},
\end{align*}
so constraint $\left(\ref{eq:5}\right)$ holds. Moreover,
\begin{align*}
    \lim_{n\to \infty}\sum_{l=1}^{\lfloor n/p \rfloor}\varepsilon_{lp,n}' & = -\lim_{n\to \infty} \sum_{l=1}^{\lfloor n/p \rfloor}\varepsilon_{lp,n}\\
    & = 0,
\end{align*}
where the latest limit is due to the assumption that $\mathbb{P}_{n}$ satisfies constraint $\left(\ref{eq:6}\right)$.
\end{proof}

\section{A Central Limit Theorem for E-K  Distributions on $\mathbb{N}$}
We would like to formalize the notion of asymptotic normality for $\omega(X_{j})$ in the limit of a sequence of random variables $\left(X_{j}\right)_{j\geq 1}$ on $\mathbb{N}$. Towards that goal we provide the following definition.

\begin{defn}\label{def:2}
Given an infinite sequence of random variables $X_1, X_2, \ldots$
defined on $\mathbb{N}$, let $X_{j}\left(n\right)$ denote the truncation of $X_{j}$ on $\left[n\right]$. Define $\varepsilon_{l,j,n}\coloneqq \mathbb{P}\left(X_{j}\left(n\right)=l\right)-1/n$. 
The sequence $X_{1}, X_{2}, \ldots$ is said to satisfy the \textbf{uniformity along primes} property if the following holds.
There exists constants $C$ and $D$ such that for all $n$ there exists a $d\geq 1$ such that for all $j\geq d$:
\begin{itemize}
\item For each prime $p>\alpha_n$
\[
\sum_{l=1}^{\left\lfloor n/p \right\rfloor }\varepsilon_{lp,j,n}\le\frac{C}{p}.
\]
\item For each $k\geq 1$ and for all $k$-tuples $\left(p_{1},\ldots,p_{k}\right)$
consisting of distinct primes of size at most $\alpha_n$
\[
\sum_{l=1}^{\left\lfloor \frac{n}{p_{1}\cdots p_{k}}\right\rfloor }\varepsilon_{lp_{1}\cdots p_{k},j,n}\le\frac{D}{n}.
\]
\end{itemize}
In addition, for any prime $p$,
\[
\lim_{j \to \infty} \lim_{n \to \infty} \sum_{l=1}^{\left\lfloor n/p \right\rfloor} \varepsilon_{lp, j, n} = 0.
\]
\end{defn}

\begin{thm}\label{thm:5}
Assume $X_{1}, X_{2}, \ldots$ is an infinite sequence of random variables defined on $\mathbb{N}$ satisfying the \textbf{uniformity along primes} property. Let $\mathbb{P}_{j,n}$ be the probability distribution of $X_{j}\left(n\right)$ on $\left[n\right]$. If $\left(n_{j}\right)_{j\in \mathbb{N}}$  is any sequence in $\mathbb{N}\setminus\lbrace 1\rbrace$ tending to $\infty$ such that $(\ref{eq:4} - \ref{eq:5})$ hold for $X_{j}\left(n_{j}\right)$ for all $j\geq 1$, then
\[
\mathbb{P}_{j,n_{j}}\left(m\le n_{j}:\omega\left(m\right)-\log\log n_{j} \le x\left(\log\log n_{j} \right)^{1/2}\right)
\to\mathbb{P}\left(Z\le x\right)
\]
as $j\to \infty$.
\end{thm}

\begin{proof}
Let's proceed by contradiction. That is, suppose $\left(n_{j}\right)_{j\in\mathbb{N}}$ is any sequence in $\mathbb{N} \setminus \left\{ 
1 \right\}$ tending to $\infty$ and suppose there exist a sequence $\left(j_{k}\right)_{k\in\mathbb{N}}$, some $x_0\in\mathbb{R}$, and some $\varepsilon_{0}>0$ such that 
\[\left| \mathbb{P}_{j_{k}, n_{j_{k}}}\left(m\leq n_{j_{k}}:\omega\left(m\right)-\log\log n_{j_{k}} \le x_{0}\left(\log\log n_{j_{k}} \right)^{1/2}\right)
- \mathbb{P}\left(Z\le x_{0}\right) \right| \geq \varepsilon_{0}
\]
for all $k\ge 1$. Consider the following PMF on $\left[n\right]$ for all $n > 1$, which satisfies the hypotheses of Theorem \ref{thm:2}:
\[ 
\mathbb{P}_{n}^{*}\left(i\right) \coloneq \begin{cases}
1/n & \text{if }n_{j_{k}}\not=n\text{ for all }k\ge1,\\
\mathbb{P}_{j_{k'\left(n\right)}, n}\left(i\right) & \text{if 
 }n_{j_{k'\left(n\right)}} = n,
\end{cases}
\]
where $k'\left( n \right)$ is the smallest integer satisfying $n_{j_{k'\left(n\right)}} = n$. It is clear that $\mathbb{P}_{n}^{*}$ is defined for all $n > 1$ and satisfies the hypotheses of Theorem \ref{thm:2} due to Definition \ref{def:2}. Let $\varepsilon = \varepsilon_{0}$, then there is some $d \ge 1$ such that
\[
\left| \mathbb{P}_{j_{k'\left(n\right)}, n_{j_{k'\left(n\right)}}}\left(m\leq n_{j_{k'\left(n\right)}}:\omega\left(m\right)-\log\log n_{j_{k'\left(n\right)}} \le x_0\left(\log\log n_{j_{k'\left(n\right)}} \right)^{1/2}\right)
- \mathbb{P}\left(Z\le x_0\right) \right| < \varepsilon
\]
if $n\geq d$, by Theorem \ref{thm:2}. The fact that at least one $k'\left(n\right)$ exists such that $n_{j_{k'\left(n\right)}} = n$ when $n\geq d$ is due to the fact that $n_{j_{k}}\to \infty$ as $k\to \infty$, so there are infinitely many $k'\left(n\right)$ such that $n_{j_{k'\left(n\right)}} = n$ when $n\geq d$. Thus, Theorem \ref{thm:5} holds by contradiction.
\end{proof}





\subsection{Moment Generating Functions and Characteristic Functions}
\begin{defn}
Let $X$ denote a random variable. The  \textbf{moment generating function (MGF)} of $X$ is defined as $M_{X}\left(t\right)=\mathbb{E}\left(e^{tX}\right)$, and the \textbf{characteristic function (CF)} of $X$ is defined as $\phi_{X} \left(t\right)=\mathbb{E}\left(e^{itX}\right)$.
\end{defn}

We have the following two facts about the MGF and the CF that will be used in the proof of Corollary \ref{cor:3}:

\begin{itemize}
\item The $k^{th}$ derivative of $M_{X}\left(t\right)$ at $t=0$ is the $k^{th}$ moment of $X$.

\item L\`evy's Continuity Theorem: Convergence in distribution, $X_{n} \overset{d}{\to} X$, for a sequence of random variables is equivalent to pointwise convergence, $\phi_{X_{n}}\to \phi_{X}$, of the corresponding CFs on all of $\mathbb{R}$.
\end{itemize}

Corollary \ref{cor:3} and Corollary \ref{cor:4} below allow us to state Theorem 7.1 of \cite{key-2} in a slightly different way.

\begin{cor}\label{cor:3}
Suppose $X_{1}, X_{2}, \ldots$ is an infinite sequence of random variables defined on $\mathbb{N}$ that satisfies the \textbf{uniformity along primes} property. Let $\left(n_{j}\right)_{j \in \mathbb{N}}$ be any sequence in $\mathbb{N} \setminus \left\{ 
1 \right\}$ tending to $\infty$ such that (\ref{eq:4} - \ref{eq:5}) hold for $X_{j}\left(n_{j}\right)$ for all $j\geq 1$. Let $\mathbb{P}_{j}$ be the probability distribution of $X_{j}$ on $\mathbb{N}$, $\mathbb{P}_{j,n}$ be the probability distribution of $X_{j}\left(n\right)$ on $\left[n\right]$, and define $\mu_{j,n} \coloneqq \mathbb{E}_{j, n}\left(\omega\left(X_{j}\left(n\right)\right)\right)$, then the following is true:
\[
\mathbb{P}_{j,n_{j}}\left(m\leq n_{j}:\omega\left(m\right)-\mu_{j,n_{j}} \le x\left(\mu_{j,n_{j}} \right)^{1/2}\right) \to\mathbb{P}\left(Z\le x\right)
\]
as $j\to \infty$.
Additionally, if $\mu_{j}<\infty$ for all $j\geq 1$ and
\[
\lim_{j\to \infty}\left(\mu_{j} / \left(\log\log n_{j}\right)^{1/2} - \left(\log\log n_{j} \right)^{1/2} \right) = 0,
\]
then
\[
\mathbb{P}_{j}\left(m:\omega\left(m\right)-\mu_{j} \le x\left(\mu_{j} \right)^{1/2}\right) \to\mathbb{P}\left(Z\le x\right)
\]
as $j\to \infty$, where
\[
\mu_{j} \coloneqq\mathbb{E}_{j}\left(\omega\left(X_{j}\right)\right) =\lim_{n\to\infty}\mu_{j,n}.
\]
We also have the following asymptotic properties:
\[
\lim_{j\to \infty}\left(\mu_{j, n_{j}} / \left(\log\log n_{j}\right)^{1/2} - \left(\log\log n_{j} \right)^{1/2} \right) = 0,
\]
\[
\lim_{j\to \infty}\left(\mu_{j, n_{j}}^{1/2} - \left(\log\log n_{j}\right)^{1/2} \right) = 0,
\]
and
\[
\lim_{j\to\infty} \frac{\mu_{j,n_{j}}}{\log \log n_{j}} = \lim_{j\to\infty} \frac{\mu_{j}}{\log \log n_{j}} = 1.
\]
\end{cor}

\begin{proof}
Let $M_{j, n_{j}}$ and $\phi_{j, n_{j}}$ be the MGF and CF of
\[
\frac{\omega\left( X_{j}\left(n_{j}
\right)\right) - \log \log n_{j}}{\sqrt{\log \log n_{j}}},
\]
and let $M_{j, \mu}$ and $\phi_{j, \mu}$ be the MGF and CF of
\[
\frac{\omega\left( X_{j}\left(n_{j}\right)\right) - \mu_{j, n_{j}}}{\sqrt{\mu_{j, n_{j}}}}.
\]
We have
\[
M_{j,n_{j}}\left(t\right) = e^{-t \left(\log\log n_{j}\right)^{1/2}} M_{\omega\left(X_{j}\left(n_{j}\right)\right)}\left(\left(\log\log n_{j}\right)^{ - 1/2} t \right)
\]
and 
\[
M_{j, \mu}\left(t\right) = e^{- t \mu_{j,n_{j}}^{1/2} } M_{\omega\left(X_{j}\left(n_{j}\right)\right)}\left(\left(\mu_{j,n_{j}}\right)^{-1/2}t\right).
\]
Let's compute the mean from the MGF $M_{j,n_{j}}\left(t\right)$:
\begin{align*}
M_{j,n_{j}}'(0) & = - \left(\log\log n_{j}\right)^{1/2}e^{-t\left(\log\log n_{j}\right)^{1/2}} M_{\omega\left(X_{j}\left(n_{j}\right)\right)}\left(\left(\log\log n_{j}\right)^{ - 1/2} t \right)\\
& \phantom{=} + e^{-t \left(\log\log n_{j}\right)^{1/2}} M^{'}_{\omega\left(X_{j}\left(n_{j}\right)\right)}\left(\left(\log\log n_{j}\right)^{ - 1/2} t \right) \left(\log\log n_{j}\right)^{-1/2} \Big\vert_{t=0}\\
& = -\left(\log\log n_{j}\right)^{1/2} +  \mu_{j,n_{j}} \left(\log\log n_{j}\right)^{-1/2}\\
& \to 0,
\end{align*}
where the last limit holds due to the proof of Theorem \ref{thm:2} (we showed the moments all approach the moments of the normal distribution) which is used in the proof of Theorem \ref{thm:5} above; thus,
\begin{equation}\label{eq:11}
\lim_{j\to \infty}\left(\mu_{j, n_{j}}/\left(\log\log n_{j}\right)^{1/2} - \left(\log\log n_{j}\right)^{1/2} \right) = 0
\end{equation}
holds.
Therefore,
\[
\lim_{j\to \infty}\frac{\mu_{j,n_{j}}}{\log\log n_{j}} = 1.
\]
From equation $\left(\ref{eq:11}\right)$ we can conclude
\[
\lim_{j\to \infty} \left( \mu_{j, n_{j}}^{1/2} - \left(\log \log n_{j}\right)^{1/2}\right) = 0.
\]
This is because, 
\[
\left| \mu_{j, n_{j}}^{1/2} - \left(\log \log n_{j}\right)^{1/2} \right| =
\left| \frac{\mu_{j, n_{j}} - \left(\log \log n_{j}\right)}{\mu_{j, n_{j}}^{1/2} + \left(\log\log n_{j} \right)^{1/2}} \right|
\leq \left| \frac{\mu_{j, n_{j}} - \log\log n_{j}}{\left(\log\log n_{j}\right)^{1/2}}\right| \to 0.
\]

Now we want to show convergence in distribution using the characteristic functions. We use the fact that the characteristic function $\phi_{\omega\left(X_{j}\left(n_{j}\right)\right)}\left(t\right)$ is continuous at $t=0$. According to properties of continuous functions, we have the ratios 
\begin{equation}\label{eq:12}
\frac{
e^{-it \mu_{j, n_{j}}^{1/2}}}
{e^{-it \left(\log \log n_{j}\right)^{1/2}}} =
e^{-it\left( \mu_{j, n_{j}}^{1/2} - \left(\log\log n_{j}\right)^{1/2}\right)} \to e^{0} = 1,
\end{equation}
and also
\begin{equation}\label{eq:13}
\frac{\phi_{\omega\left(X_{j}\left(n_{j}\right)\right)}\left(\left(\mu_{j, n_{j}}\right)^{-1/2}t\right)}{\phi_{\omega\left(X_{j}\left(n_{j}\right)\right)}\left(\left(\log \log n_{j} \right)^{ - 1/2} t \right)
} \to \frac{\phi_{\omega\left(X_{j}\left(n_{j}\right)\right)}\left( 0 \right)}{\phi_{\omega\left(X_{j}\left(n_{j}\right)\right)}\left( 0 \right)} = 1
\end{equation}
as $j\to \infty$. This gives us the required convergence of $\phi_{j, \mu}$ at any $t$; i.e.,
\begin{align*}
\phi_{j, \mu}\left( t \right) = e^{-it \mu_{j, n_{j}}^{1/2}}
\phi_{\omega\left(X_{j}\left(n_{j}\right)\right)}\left(\left(\mu_{j, n_{j}}\right)^{-1/2}t\right)
& \overset{\left(\ref{eq:12},\ref{eq:13}\right)}\sim e^{-it \left(\log\log n_{j}\right)^{1/2}} \phi_{\omega\left(X_{j}\left(n_{j}\right)\right)}\left(\left(\log\log n_{j}\right)^{ - 1/2} t \right)\\
& \to e^{-t^2/2} = \phi_Z\left(t \right).
\end{align*}
By L\`evy's Continuity Theorem,
\[
\mathbb{P}_{j,n_{j}}\left(m\leq n_{j}:\omega\left(m\right)-\mu_{j,n_{j}} \le x\left(\mu_{j,n_{j}} \right)^{1/2}\right) \to\mathbb{P}\left(Z\le x\right).
\]

Suppose additionally that $\mu_{j}<\infty$ for all $j\geq 1$ and
\[
\lim_{j\to \infty}\left(\mu_{j} / \left(\log\log n_{j}\right)^{1/2} - \left(\log\log n_{j} \right)^{1/2} \right) = 0.
\]
Similarly, let $\phi_{j,\mu_{j}}$ be the characteristic function of
\[
\frac{\omega\left( X_{j}\right) - \mu_{j}}{\sqrt{\mu_{j}}}.
\]
Then
\begin{align*}
\phi_{j, \mu_{j}}\left( t \right) = e^{-it\mu_{j}^{1/2}}
\phi_{\omega\left(X_{j}\right)}\left(\left(\mu_{j}\right)^{-1/2}t \right)
& \sim e^{-it \left(\log\log n_{j}\right)^{1/2}} \phi_{\omega\left(X_{j}\left(n_{j}\right)\right)}\left(\left(\log\log n_{j}\right)^{-1/2}t\right)\\
& \to e^{-t^2/2} = \phi_Z\left(t \right),
\end{align*}
for any $t\in \mathbb{R}$. Thus, again we have that
\[
\mathbb{P}_{j}\left(m:\omega\left(m\right)-\mu_{j} \le x\left(\mu_{j} \right)^{1/2}\right) \to\mathbb{P}\left(Z\le x\right),
\]
by L\`evy's Continuity Theorem. This completes the proof.
\end{proof}

\begin{rem*}
We can replace $\mu_{j}$ in Corollary \ref{cor:3} with any sequence $a_{j}$ that satisfies the limit 
\[
\lim_{j\to \infty}\left(a_{j} / \left(\log\log n_{j}\right)^{1/2} - \left(\log\log n_{j} \right)^{1/2} \right) = 0,
\]
to get a similar result. This is a generalization of Theorem 7.1 in \cite{key-2} which states this result for particular sequences of $Zeta\left(s \right)$ distributions as their parameter $s\to 1$. In their result, the mean was $\sum_{p} p^{-s}$, where the sum is over all prime numbers.
\end{rem*}

\begin{defn}\label{def:3}
Given an infinite sequence of random variables $X_{1}, X_{2}, \ldots$
defined on $\mathbb{N}$, let $X_{j}\left(n\right)$ denote the truncation of $X_{j}$ on $\left[n\right]$.
If $X_{1}, X_{2}, \ldots$ satisfies the uniformity along primes property, and $X_{j}\left(n\right)$ satisfies (\ref{eq:4} - \ref{eq:5}) for any $n$ and any $j$, then we say the sequence satisfies the \textbf{complete uniformity along primes} property.
\end{defn}

The following lemma will be applied in Corollary $\ref{cor:4}$.
\begin{lem}\label{lem:3}
If $X_1, X_2, \ldots$ has the \textbf{complete uniformity along primes} property and $\mu_{j} < \infty$ for all $j\geq 1$, then $\mu_{j}\to \infty$ as $j\to \infty$.
\end{lem}
\begin{proof}
For any sequence $\left(n_{j}\right)_{j \in \mathbb{N}}$ in $\mathbb{N} \setminus \left\{ 
1 \right\}$  tending to  $\infty$ we have $\log\log n_{j} \sim \mu_{j, n_{j}} \sim \mu_{j}$ as $j\to \infty$.
\end{proof}

\begin{cor}\label{cor:4}
If $X_{1}, X_{2}, \ldots$ has the \textbf{complete uniformity along primes} property and $\mu_{j} < \infty$ for all $j\geq 1$, then
\[
\mathbb{P}_{j}\left(m:\omega\left(m\right)-\mu_{j} \le x\left(\mu_{j}\right)^{1/2}\right)\to\mathbb{P}\left(Z\le x\right)
\]
as $j\to \infty$.
\end{cor}

\begin{proof}
Define $n_{j} \coloneqq \max\left\{ \left\lfloor e^{e^{\mu_{j}}}\right\rfloor, 2 \right\}$ so that $n_{j}\to\infty$ by Lemma \ref{lem:3} and the limit
\[
\lim_{j\to \infty}\left(\mu_{j} / \left(\log\log n_{j}\right)^{1/2} - \left(\log\log n_{j} \right)^{1/2} \right) = 0
\]
in Corollary \ref{cor:3} holds.
\end{proof}

\subsection{The Zeta Distribution}
We show the zeta distribution has the \textbf{complete uniformity along primes} property as $s\to 1$. Given $s>1$, denote by $Z_{s}$ the $Zeta\left(s\right)$ distribution
so that for any $j\in\mathbb{N}$, 
\[Z_{s}\left(j\right)=\frac{1}{j^{s}\zeta\left(s\right)},\]
where 
\[\zeta\left(s\right)=\sum_{j\ge1}\frac{1}{j^{s}}\] 
denotes the Riemann zeta function. Since Theorem \ref{thm:2} involves distributions defined
on $\left[n\right]$, restrict the $Zeta\left(s\right)$ distribution
to $\left[n\right]$ and then normalize by dividing by $\sum_{i=1}^{n}\frac{1}{i^{s}\zeta\left(s\right)}$;
i.e., for $j\in\left[n\right]$,
\begin{align*}
Z_{s,n}\left(j\right) & \coloneqq\frac{\frac{1}{j^{s}\zeta\left(s\right)}}{\sum_{i=1}^{n}\frac{1}{i^{s}\zeta\left(s\right)}}\\
 & =\frac{1}{j^{s}\sum_{i=1}^{n}\frac{1}{i^{s}}};
\end{align*}
and $Z_{s,n}$ is known as the Zipf distribution with parameters $n$
and $s$. For a Zifp distribution, we have 
\[ \varepsilon_{i,n} \overset{\left(\ref{eq:1}\right)}{=} \frac{1/i^{s}}{\sum_{j=1}^{n}1/j^s} - \frac{1}{n}.\]
Thus,
\begin{align*}
\sum_{l=1}^{\lfloor n/p \rfloor}\varepsilon_{lp,n} & =\frac{\sum_{l=1}^{\lfloor n/p\rfloor}1/\left(lp\right)^{s}}{\sum_{l=1}^{n}1/l^{s}}-\frac{\lfloor n/p\rfloor}{n}\\
 & =\frac{1}{p^{s}}\frac{\sum_{l=1}^{\lfloor n/p\rfloor}1/l^{s}}{\sum_{l=1}^{n}1/l^{s}}-\frac{\lfloor n/p\rfloor}{n}\\
 & \le \frac{1}{p^s}\\
 & \le \frac{1}{p},
\end{align*}
so constraint $\left(\ref{eq:4}\right)$ holds with $C=1$ for any $s$. Moreover,
\begin{align*}
\sum_{l=1}^{\left\lfloor \frac{n}{p_{1}\cdots p_{k}}\right\rfloor }\varepsilon_{lp_{1}\cdots p_{k},n} & =\frac{\sum_{l=1}^{\lfloor \frac{n}{p_{1}\cdots p_{k}}\rfloor}1/\left(lp_{1}\cdots p_{k}\right)^{s}}{\sum_{l=1}^{n}1/l^{s}}-\frac{\lfloor \frac{n}{p_{1}\cdots p_{k}}\rfloor}{n}\\
 & =\frac{1}{\left(p_{1}\cdots p_{k}\right)^{s}}\frac{\sum_{l=1}^{\lfloor n/p\rfloor}1/l^{s}}{\sum_{l=1}^{n}1/l^{s}}-\frac{1}{p_{1}\cdots  p_{k}}+\frac{1}{n}\\
 & \le \frac{1}{\left(p_{1}\cdots p_{k}\right)^{s}} - \frac{1}{p_{1}\cdots p_{k}}+\frac{1}{n}\\
 & \le \frac{1}{n},
\end{align*}
so constraint $\left(\ref{eq:5}\right)$ holds with $D=1$ for any $s$. Furthermore,
\begin{align*}
\lim_{n\to \infty}\sum_{l=1}^{\lfloor n/p \rfloor}\varepsilon_{lp,n} & =\lim_{n\to \infty}\frac{\sum_{l=1}^{\lfloor n/p\rfloor}1/\left(lp\right)^{s}}{\sum_{l=1}^{n}1/l^{s}}-\frac{1}{p}\\
 & =\frac{1}{p^{s}}-\frac{1}{p}\\
& \overset{s\to 1}{\to} 0;
\end{align*}
therefore, $\left(\ref{eq:6}\right)$ holds as $s\to 1$.

\begin{cor}
Let $\left(a_{j}\right)_{j \in \mathbb{N}}$ be any sequence of real numbers such that $0 < a_{j}\to \infty$. If $X_{j}$ is a sequence of $Zeta\left(1+1/a_{j}\right)$-distributed random variables, then 
\[\mu_{j} = \sum_{p} p^{-\left(1+1/a_{j}\right)}<\infty,\]  
\[\lim_{j \to \infty}\mu_{j}=\infty,\] and
\[
\mathbb{P}_{j}\left(m:\omega\left(m\right)-\mu_{j} \le x\left(\mu_{j} \right)^{1/2}\right) \to\mathbb{P}\left(Z\le x\right)
\]
as $j\to \infty$.
\end{cor}
This latest limit is the statement of Theorem 7.1 in \cite{key-2}.

\subsection{Convex Combinations of E-K Distributions on $\mathbb{N}$}

We can also take convex combinations of sequences to form new sequences which satisfy the uniformity along primes property.

\begin{cor}\label{cor:6}
Let $\lambda \in \left[0,1\right]$. Let $X_{1}, X_{2},\ldots$ and $Y_{1}, Y_{2},\ldots$ be two sequences of random variables on $\mathbb{N}$ that both satisfy the \textbf{complete uniformity along primes} property. Then we can define a new sequence $
Z_{1}, Z_{2}, \ldots$
such that each $Z_{j}$ is the convex combination of $X_{j}$ and $Y_{j}$. In particular, 
\[
\mathbb{P}_{Z_{j}}\left(i\right) \coloneqq \lambda \mathbb{P}_{X_{j}}\left(i\right) + \left(1-\lambda\right)\mathbb{P}_{Y_{j}}\left(i\right); i\geq 1, j\geq 1.
\]

Then $Z_{j}$ satisfies the \textbf{complete uniformity along primes} property and thus Theorem \ref{thm:5} holds for the $Z_{j}$ sequence.
\end{cor}

\begin{proof}
The proof is similar to that of Theorem \ref{thm:3}.
\end{proof}

\begin{cor}
Apply Corollary \ref{cor:6} to any convex combination of a finite number of sequences
\[
\left(X_{1,1}, X_{1,2}, \ldots \right), \ldots, \left(X_{n,1}, X_{n,2}, \ldots\right)
\]
that all satisfy the \textbf{complete uniformity along primes} property. Then the conclusion of Corollary \ref{cor:6} holds for this convex combination.
\end{cor}

\begin{proof}
Use induction and Corollary \ref{cor:6}.
\end{proof}

\section{Illustrative Examples}
In this section, we give examples from the class of distributions $\mathbb{P}^{*}_n$ that satisfy the hypotheses $\left(\ref{eq:4} - \ref{eq:6}\right)$ and examples of limits of distributions to which we can apply Theorem \ref{thm:5}. It will be shown that the statement of Theorem \ref{thm:2} holds when $\mathbb{P}_{n}^{*}$ is replaced with either the Harmonic$\left(n\right)$ distribution or a convex combination of Harmonic and Uniform Distributions. Then we show that the $Zeta\left(s\right)$ and the Logarithmic distribution satisfy Theorem \ref{thm:5} as their parameters tend towards limits. We introduce a $2$-parameter family of distributions $LZ\left(s,\alpha\right)$ which includes the Logarithmic distribution and Zeta distribution as special cases (when $\alpha=1$ and $s = 1$ respectfully). We also look at a geometric power series distribution that converges to the normal distribution on all truncations as $s \to 1$.

\subsection{The Harmonic Distribution}
This was proved in $\S\S$2.1 above.

\subsection{The Zeta Distribution}
This was proved in $\S\S$4.2 above.

\subsection{The Logarithmic Distribution} Given a real number $s$ with $0<s<1$, a logarithmic distribution with parameter $s$ is given by 
\[L_{s}\left(i\right) \coloneq \frac{-1}{\log \left(1-s\right)} \frac{s^i}{i}; \,\, i \in \mathbb{N}.\]
Now we will show that the number of distinct prime factors, $\omega\left(\cdot\right)$,
of a truncated log-distributed variable has the same central limit theorem as the uniform distribution as $n\to \infty$ and $s\to 1$. For $i\in\left[n\right]$ we have 
\begin{align*}
L_{s,n}\left(i\right) & \coloneq \frac{\frac{-1}{\log\left(1-s\right)}\frac{s^{i}}{i}}{\sum_{l=1}^{n}\frac{-1}{\log\left(1-s\right)}\frac{s^{l}}{l}}=\frac{\frac{s^{i}}{i}}{\sum_{l=1}^{n}\frac{s^{l}}{l}},
\end{align*}
so 
\begin{align*}
\varepsilon_{i,n} & 
\overset{\left(\ref{eq:1}\right)}{=}
\frac{\frac{s^{i}}{i}}{\sum_{l=1}^{n}\frac{s^{l}}{l}}-\frac{1}{n}.
\end{align*}
Therefore,
\begin{align*}
\sum_{l=1}^{\left\lfloor \frac{n}{p}\right\rfloor }\varepsilon_{lp,n} & =\sum_{l=1}^{\left\lfloor \frac{n}{p}\right\rfloor }\left(\frac{\frac{\left(s\right)^{lp}}{lp}}{\sum_{l=1}^{n}\frac{s^{l}}{l}}-\frac{1}{n}\right)\\
 & \le\frac{1}{p}\frac{\sum_{l=1}^{\left\lfloor \frac{n}{p}\right\rfloor }\frac{\left(s^{p}\right)^{l}}{l}}{\sum_{l=1}^{n}\frac{s^{l}}{l}}-\frac{\lfloor n/p\rfloor}{n}\\
 & \le\frac{1}{p}-\frac{\lfloor n/p\rfloor}{n}\\
 & \le\frac{1}{p},
\end{align*}
so $\left(\ref{eq:4}\right)$ holds with $C=1$. Similarly, 
\begin{align*}
\sum_{l=1}^{\left\lfloor \frac{n}{p_{1}\cdots p_{k}}\right\rfloor }\varepsilon_{lp_{1}\cdots p_{k},n} & =\sum_{l=1}^{\left\lfloor \frac{n}{p_{1}\cdots p_{k}}\right\rfloor }\left(\frac{\frac{\left(s\right)^{lp_{1}\cdots p_{k}}}{lp_{1}\cdots p_{k}}}{\sum_{l=1}^{n}\frac{s^{l}}{l}}-\frac{1}{n}\right)\\
 & \le\frac{1}{p_{1}\cdots p_{k}}\frac{\sum_{l=1}^{\left\lfloor \frac{n}{p_{1}\cdots p_{k}}\right\rfloor }\frac{\left(s^{p_{1}\cdots p_{k}}\right)^{l}}{l}}{\sum_{l=1}^{n}\frac{s^{l}}{l}}-\frac{1}{p_{1}\cdots p_{k}}+\frac{1}{n}\\
 & \le\frac{1}{p_{1}\cdots p_{k}}-\frac{1}{p_{1}\cdots p_{k}}+\frac{1}{n}\\
 & =\frac{1}{n},
\end{align*}
so $\left(\ref{eq:5}\right)$ holds with $D=1$. Moreover, we have
\[
\sum_{l=1}^{\lfloor n/p\rfloor}\varepsilon_{lp,n} \overset{n\to \infty}{\to} \frac{1}{p}\frac{\log\left(1-s^p \right)}{\log\left(1-s\right)} - \frac{1}{p} \overset{s\to 1}{\to} 0.
\]
Therefore, $\left(\ref{eq:6}\right)$ holds, so $\omega\left(L_{i,s}\right)$ is asymptotically distributed as $\mathcal{N}\left(\log \log n, \log \log n\right)$ as $n\to \infty$ and $s\to 1$.

Although Kac's heuristic for Theorem \ref{thm:1} is based on the asymptotic independence in the uniform case, we will show that the events $\lbrace \text{divisible by } p\rbrace$ and $\lbrace \text{divisible by }q\rbrace$ are not independent for any value of $s$. Let $A_{p}$ denote the set of all positive integers divisible by
$p$, then 
\begin{align*}
L_{s}\left(A_p\right) & = \sum_{l=1}^{\infty} \frac{-1}{\log \left(1-s\right)} \frac{\left(s^p\right)^l}{pl} \\
 & =\frac{1}{p}\frac{\log\left(1-s^{p}\right)}{\log\left(1-s\right)}
\end{align*}
and 
\begin{align*}
L_{s}\left(A_{p}\cap A_{q}\right) & = \sum_{l=1}^{\infty} \frac{-1}{\log \left(1-s\right)} \frac{\left(s^{pq}\right)^l}{pql} \\
 & =\frac{1}{pq}\frac{\log\left(1-s^{pq}\right)}{\log\left(1-s\right)}.
\end{align*}
It is worth noting that 
\[ \lim_{s\to 1}L_{s}\left(A_{p} \cap A_{q}\right) =  \frac{1}{pq} = \lim_{s\to 1}\left(L_{s}\left(A_{p}\right) L_{s}\left(A_{q}\right)\right),\] as $s\to 1$; so independence is only approached in the limit.

\subsection{Geometric Power Series Distribution}
Let $s\in\left(0,1\right)$. Define:
\[
\mathbb{P}\left(i\right) \coloneqq \dfrac{1-s}{s}s^{i}; \,\, i \in \mathbb{N}
\]
then this distribution satisfies the hypotheses of Theorem 2 in the limit as $s\to 1$. As $s\to 1$ it seems that this distribution is converging to the uniform distribution on any truncation, so the result is to be expected.

Truncating on $\left[n\right]$ leads to
\[
\varepsilon_{i,n}
\overset{\left(\ref{eq:1}\right)}{=}\frac{s^{i}}{\sum_{j=1}^{n} s^j}-\frac{1}{n},
\] 
therefore
\begin{align*}
\sum_{l=1}^{\lfloor  n/p \rfloor} \varepsilon_{lp,n} & = \sum_{l=1}^{\lfloor n/p \rfloor}\left(\frac{s^{i}}{\sum_{j=1}^{n}s^{j}} -\frac{1}{n}\right)\\
& = \frac{s-s^{\lfloor n/p\rfloor +1}}{s-s^{n+1}}-\frac{\lfloor n/p \rfloor}{n}\\
& \overset{s\to 1}{\to} \frac{\lfloor n/p \rfloor }{n} - \frac{\lfloor n/p \rfloor}{n}\\
& =0,
\end{align*}
so constraints $\left(\ref{eq:4}-\ref{eq:6}\right)$ hold with $C=D=1$. 

\subsection{Convex Combination of Harmonic and Uniform}

Following Cranston and Mountford in \cite{key-1}, let $\lambda \in [0,1]$ and define 
\[
\mathbb{P}_n(i) \coloneqq \frac{\lambda}{i h_{n}} + \left(1-\lambda\right)\frac{1}{n}; \,\, i \in \mathbb{N}
\]
where $h_{n}$ is the $n^{\text{th}}$ harmonic number. Then $\mathbb{P}_n$ satisfies the hypotheses of Theorem \ref{thm:2}. This is proved in Corollary \ref{cor:2} above using Theorems \ref{thm:2} and \ref{thm:3}.

\subsection{A Logarithmic-Zeta Distribution}

There is a $2$-parameter family of power series distributions on $\mathbb{N}$ for $(s,\alpha)\in\mathbb{R}^2$ such that $0<s\leq1$, $\alpha\geq 1$, and $s=\alpha=1$ is not allowed; the PMF is given by
\[LZ_{s,\alpha}\left(i\right) \coloneq \frac{1}{\sum_{j=1}^{\infty} \frac{s^j}{j^\alpha}} \frac{s^i}{i^{\alpha}}; \,\, i \in \mathbb{N}.\]
Similarly, there is a truncation of $LZ$ given by
\[LZ_{s,\alpha,n}\left(i\right) \coloneq \frac{1}{\sum_{j=1}^{n} \frac{s^j}{j^\alpha}} \frac{s^i}{i^{\alpha}}; \,\, i \leq n,\]
and a description of the $\varepsilon_{i,n}$ given by
\[
\varepsilon_{i,n} \overset{\left(\ref{eq:1}\right)}{=} \frac{1}{\sum_{j=1}^{n} \frac{s^j}{j^\alpha}} \frac{s^i}{i^{\alpha}} - \frac{1}{n}.
\]
It is clear that a similar type of argument shows that Theorem \ref{thm:2} holds for the truncated logarithmic-zeta distribution as $(s, \alpha) \to (1,1)$; in particular, we obtain 
\[
\sum_{l=1}^{\lfloor n/p \rfloor}\varepsilon_{lp,n} \le \frac{1}{p^{\alpha}} \leq \frac{1}{p}
\]
and
\[
\sum_{l=1}^{\lfloor \frac{n}{p_{1}\cdots p_{k}} \rfloor}\varepsilon_{lp_{1}\cdots p_{k},n} \le \frac{1}{p_{1}^{\alpha} \cdots p_{k}^{\alpha}}-\frac{1}{p_{1} \cdots p_{k}}+\frac{1}{n}\leq\frac{1}{n}.  
\]
Furthermore, 
\[\lim_{n\to \infty}
\sum_{l=1}^{\lfloor n/p \rfloor}\varepsilon_{lp,n} \overset{s \to 1}{\to} 
\frac{1}{p^{\alpha}} - \frac{1}{p}
\overset{\alpha \to 1}{\to} 0.
\]

\subsection{Passage to the Limit: 
$\lim\limits_{\substack{\mathllap{s} \to \mathrlap{1} \\ \mathllap{\alpha} \to \mathrlap{1}}}LZ_{s,\alpha}$}
We would like to look at the behavior of $\omega \left( X \right)$ for $X$ distributed as $LZ_{s,\alpha}$ as $\left(s,\alpha\right) \to \left(1,1\right)$. The motivation for this comes from the fact that as $(s,\alpha) \to (1,1)$ the truncated distributions resemble the harmonic distribution which we know behaves similarly to the uniform distribution on $\left[n\right]$ when $n\to\infty$ according to Theorem \ref{thm:2}.

In Theorem 5.5 of \cite{key-2}, they provide the moment generating function for $\omega\left(X_{\alpha}\right)$ when $X_{\alpha}$ is a $Zeta\left(\alpha\right)$ distributed random integer with parameter $\alpha >1$. They go on to prove that
\[
\hat{\omega}\left(X_\alpha \right) = 
\frac{\omega\left(X_{\alpha}\right) - \sum_p p^{-\alpha } }{\sqrt{\sum_p p^{-\alpha }}} \overset{d}{\to} Z
\]
as $\alpha \to 1$ in Theorem 7.1 of \cite{key-2} by using this moment generating function.
Later, Cranston and Mountford \cite{key-1} give a new proof of Theorem \ref{thm:1}. The proof uses Theorem 7.1 about zeta distributions to prove Erd\H{o}s-Kac in a way that translates over to settings where zeta functions still make sense.

We generalize Peltzer and Cranston's Theorem 7.1 in \cite{key-2} in the following way. As $(s,\alpha)\to \left(1,1\right)$, we show that
\[
\hat{\omega}\left(X_{s, \alpha} \right) =
\frac{\omega\left(X_{s, \alpha} \right) - \mu_{s,\alpha}}{\sqrt{\mu_{s,\alpha}}} \overset{d}{\to} Z,
\]
where $\mu_{s,\alpha}$ is the mean of $\omega\left(X_{s,\alpha}\right)$ when $X_{s,\alpha}$ is distributed as $LZ_{s,\alpha}$. In fact, we do not need independence, nor do we need to compute the mean for any particular $s,\alpha$ in order to conclude this (neither do we need to compute a MGF).

In $\S$4, we proved a stronger statement than the above statement about $LZ_{s,\alpha}$; we proved that in general
\[
\hat{\omega}\left(X_{j}\right) =
\frac{\omega\left(X_{j}\right) - \mu_{j}}{\sqrt{\mu_{j}}} 
\overset{d}{\to} Z,
\]
as long as the truncated variables $X_{j}\left(n\right)$ have the \textbf{complete uniformity along primes} property. Here we assume that $\mu_{j}=\mathbb{E}(\omega\left(X_{j}\right))<\infty$. When $\left(s,\alpha\right)\to \left(1,1\right)$, we recover the above statement about $LZ_{s, \alpha}$.

\subsection{A Non-Example: Zeroing at a Set of Primes} Fix $n\in \mathbb{N}$ and let $p\le n$ denote a prime. Consider the PMF defined by
\[ \mathbb{P}_{n,p}\left(i\right)=
\begin{cases} 
      \frac{1}{\#\left(\left[n\right] \setminus p\mathbb{N}\right)} & i \not \in p\mathbb{N}, \\
      0 & i \in p\mathbb{N}.
   \end{cases}
\]
We have
\[\varepsilon_{i,n} = \begin{cases} 
     \frac{1}{\#\left(\left[n\right] \setminus p\mathbb{N}\right)}-1/n & i \not \in p\mathbb{N}, \\
      -1/n & i \in p\mathbb{N}.
   \end{cases}
\]

\noindent Therefore,
\begin{align*}
\lim_{n\to\infty} \sum_{l=1}^{\lfloor n/p \rfloor} \varepsilon_{lp, n} & = \lim_{n \to \infty}-\frac{\lfloor n/p \rfloor}{n}\\
& = -1/p\\
& \neq 0.
\end{align*}
Thus, this PMF does not satisfy constraint $\left(\ref{eq:6}\right)$. We conjecture that the conclusion of Theorem \ref{thm:2} does not hold for this distribution.


\subsection{An Erd\H{o}s-Kac Theorem for Continuous Variables}
Consider a continuous uniform random variable $N_n$ on the interval $\left(0,n\right]$. Then $\left\lceil N_{n}\right\rceil$ is a uniform variable on $\left[n\right]$. Therefore, we have
\begin{cor}
Let $Z$ denote a standard normal variable, and let $X_{n}$ be a continuous random variable on $\left(0,n\right]$. 
Define the $\varepsilon_{i,n}$ according to the following relation: 
$\mathbb{P}_{n}\left( \left\lceil X_{n}\right\rceil=i \right)=\frac{1}{n}+\varepsilon_{i,n}$.
If the constraints
\begin{itemize}
\item There exists a constant $C$ such that for all $n > 1$ and for all primes
$p$ with $p>\alpha_n$,
\[
\sum_{l=1}^{\left\lfloor n/p \right\rfloor }\varepsilon_{lp,n}\le\frac{C}{p}.
\]
\item There exists a constant $D$ such that 
\[
\sum_{l=1}^{\left\lfloor \frac{n}{p_{1}\cdots p_{k}}\right\rfloor }\varepsilon_{lp_{1}\cdots p_{k},n}\le\frac{D}{n}
\]
for all $n > 1$ and, for each $k$, all $k$-tuples $\left(p_{1},\ldots,p_{k}\right)$
consisting of distinct primes of size at most $\alpha_n$, and 

\item For any prime $p$, 
\[
\lim_{n \to \infty} \sum_{l=1}^{\lfloor n/p \rfloor} \varepsilon_{lp, n} = 0
\]
\end{itemize}
all hold, then
\begin{align*}
\mathbb{P}_{n}\left(t\le n:\omega\left(\lceil t\rceil\right)-\log\log n\le x\left(\log\log n\right)^{1/2}\right) & \to\mathbb{P}\left(Z\le x\right)
\end{align*}
as $n\to \infty$.
\end{cor}

\section{Conclusion}
Theorem \ref{thm:2} generalizes the Erd\H{o}s-Kac Theorem for $\omega\left(\cdot \right)$ to distributions other than the uniform distribution, and this theorem was proved by imposing constraints $\left(\ref{eq:4}-\ref{eq:6}\right)$ on a PMF of the form $\mathbb{P}\left(i\right)=1/n+\varepsilon_{i,n}$. We showed that the uniform and harmonic distributions satisfy these constraints; then, we showed that any convex sum of these PMFs also satisfies the constraints (4 - 6). 

The uniformity property provides a natural way to examine asymptotic properties of truncations of variables with support $\mathbb{N}$. Given an infinite sequence $X_{1}, X_{2},\ldots$ of random variables on $\mathbb{N}$ satisfying uniformity along primes, Theorem \ref{thm:5} showed that for any sequence with $n_{j}\to \infty$ as $j\to \infty$, the distribution of $\omega\left(X_{j}\left(n_{j}\right)\right)$ is asymptotically normally distributed with mean and variance both equal to $\log \log n_{j}$ as long as $X_{j}\left(n_{j}\right) $ satisfies (\ref{eq:4} - \ref{eq:6}).

The definition of complete uniformity along primes allows us to obtain central limit theorems regardless of how $n\to\infty$, and allows us to make asymptotic statements involving the mean of $\omega\left(X_{j}\right)$. This generalized a statement from \cite{key-2} involving the mean of $\omega\left(X_s\right)$ as $s\to 1$ when $X_s$ is a random $Zeta\left(s\right)$-distributed variable.

Another way to generalize Theorem \ref{thm:2} would be to incorporate it with
other generalizations, e.g., \cite{key-1,key-2,key-7, key-8}. By incorporating Theorem \ref{thm:2} with these, further generalizations can be made in which the original setting is not $\left[n\right]$, the underlying distribution of the random-integer is not uniform, and $\omega\left(n\right)$ can be replaced with a more general strongly additive function
$f\left(n\right)$.

We also showed the complete uniformity property holds, and thus normality in the limit, for $Zeta\left(s\right)$ and a number of similar distributions.
It is suspected by the authors, but not known, whether or not the hypotheses $\left(\ref{eq:4}-\ref{eq:6}\right)$ are necessary and sufficient for the conclusion of Theorem \ref{thm:2}; we conjecture that is the case.

\end{document}